\documentclass[a4paper]{article}
\usepackage[T1]{fontenc}
\usepackage[utf8]{inputenc}
\usepackage{amsmath}
\usepackage{amsfonts}
\usepackage{amssymb}
\usepackage{graphicx}
\usepackage{tikz}
\usepackage{mathtools}
\usepackage[top=3cm, bottom=3cm, left=3cm, right=3cm]{geometry}
\usepackage{calrsfs}
\usepackage{bbm}
\usepackage{algorithm}
\usepackage{cases}
\usepackage[noend]{algpseudocode}
\usepackage{float}
\usepackage{enumerate}
\usepackage{stmaryrd}
\DeclareMathAlphabet{\pazocal}{OMS}{zplm}{m}{n}
\usepackage{lipsum}

\newcommand\blfootnote[1]{%
  \begingroup
  \renewcommand\thefootnote{}\footnote{#1}%
  \addtocounter{footnote}{-1}%
  \endgroup
}

\makeatletter
\newcommand*\bigcdot{\mathpalette\bigcdot@{.5}}
\newcommand*\bigcdot@[2]{\mathbin{\vcenter{\hbox{\scalebox{#2}{$\m@th#1\bullet$}}}}}
\makeatother

\newcommand{\mydef}{\vcentcolon=}

\newtheorem{theorem}{\indent Theorem}[section]
\newtheorem{proposition}[theorem]{\indent Proposition}
\newtheorem{assumption}[theorem]{\indent Assumption}

\newtheorem{definition-theorem}[theorem]{\indent Definition-Theorem}
\newtheorem{corollary}[theorem]{\indent Corollary}

\newenvironment{proof}{\paragraph{Proof:}}{\hfill$\square$}

\def \N{\mathbb{N}}
\def \R{\mathbb{R}}
\def \Z{\mathbb{Z}}
\def \V{\mathbb{V}}

\def \L{\mathcal{L}}

\def \P{\mathbb{P}}
\newcommand{\E}{\mathbb{E}}

\title{\bf The Effect of Graph Connectivity on Metastability in a Stochastic System of Spiking Neurons}

\author{Morgan André and Léo Planche \\ \textit{Instituto de Matemática e Estatística,} \\ \textit{Universidade de São Paulo.}}

\begin{document}

\maketitle

\blfootnote{© 2020. This manuscript version is made available under the CC-BY-NC-ND 4.0 license}

\begin{abstract}
We consider a continuous-time stochastic model of spiking neurons originally introduced by Ferrari et al. in \cite{ferrari}. In this model, we have a finite or countable number of neurons which are vertices in some graph $G$ where the edges indicate the synaptic connection between them. We focus on metastability, understood as the property for the time of extinction of the network to be asymptotically memory-less, and we prove that this model exhibits two different behaviors depending on the nature of the specific underlying graph of interaction $G$ that is chosen. In this model the spiking activity of any given neuron is represented by a point process, whose rate fluctuate between $1$ and $0$ over time depending on whether the membrane potential is positive or null. The membrane potential of each neuron evolves in time by integrating all the spikes of its adjacent neurons up to the last spike of the said neuron, so that when a neuron spikes, its membrane potential is reset to $0$ while the membrane potential of each if its adjacent neurons is increased by one unit. Moreover, each neuron is exposed to a leakage effect, modeled as an abrupt loss of membrane potential which occurs at random times driven by a Poisson process of some fixed rate $\gamma$. It was previously proven that when the graph $G$ is the infinite one-dimensional lattice, this model presents a phase transition with respect to the parameter $\gamma$. It was also proven that, when $\gamma$ is small enough, the renormalized time of extinction (the first time at which all neurons have a null membrane potential) of a finite version of the system converges in law toward an exponential random variable when the number of neurons goes to infinity. The present article is divided into two parts. First we prove that, in the finite one-dimensional lattice, this last result doesn't hold anymore if $\gamma$ is large enough,  and in fact we prove that for $\gamma > 1$ the renormalized time of extinction is asymptotically deterministic. Then we prove that conversely, if $G$ is the complete graph, the result of metastability holds for any positive $\gamma$.
\end{abstract}

\vspace{0.4 cm}



\section{Introduction}

In this paper we study the asymptotic behavior of a system of interacting point processes introduced in \cite{ferrari}, which is aimed to offer a good model for neural interactions.

\vspace{0.4 cm}

Biological neural networks consist in neurons which communicate together via synapses. For each neuron there is a difference in electrical potential between the interior of the cell and its exterior, and this difference is called the \textit{membrane potential} of the neuron. When the membrane potential of a given neuron is high enough, it tends to suddenly depolarize to a resting value. This phenomenon is known as an \textit{action potential}, or a \textit{spike}. When a spike occurs, the neurons which are linked to the neuron that just spiked, named its \textit{postsynaptic neurons}, have their membrane potential increased. The neuron that just spiked is said to be a \textit{presynaptic neuron} for the neurons it influences. Finally, between spikes, the membrane potential of the neurons tends to naturally decrease down to the resting value, a phenomenon known as \textit{leakage}.

\vspace{0.4 cm}

Roughly speaking, the model we consider consists of a graph whose nodes represents neurons and the edges the synapses between them. The membrane potential of each neuron is represented by an integer which evolves in time. The spiking activity of each neuron is represented by a continuous time point process, whose rate evolves in time with respect to the membrane potential of the neuron. Moreover the leakage effect is modeled by a sudden loss of membrane potential occurring at random times driven by a Poisson process of some parameter $\gamma$. The evolution of membrane potential of each neuron depends on the spiking activity of its presynaptic neurons up to the last event by which the said neuron has been affected (a leakage or a spike), so that the membrane potential of this neuron is also reset to $0$ whenever it spikes.

\vspace{0.4 cm}

The idea of using stochastic processes to model biological neural networks dynamics has now a long history and the reader may refer to \cite{tuckwell} for a general historical review. Inside these stochastic models, the subcategory making use of point processes has become increasingly popular during the last decades (considering, for example, Hawkes processes \cite{hawkes, chevallier, chornoboy, hansen} as well as Wold processes \cite{kass}). In the model considered here, the membrane potential is reset to $0$ after each spike, which makes it closer to the model introduced in \cite{brillinger}. Moreover it implies the biologically motivated fact that each of the spiking point processes alone has a memory of variable length. It can be seen as a continuous time version of the model introduced in \cite{galves}, sometimes named the Galves-Löcherbach model, or simply GL model. Other versions of this model has been studied since it has been introduced, see for example \cite{galves2019, papageorgiou, demasi, duarte, fournier}. We refer to \cite{review} for a general review.

\vspace{0.4 cm}

Two crucial area of research in the study of biological neural networks are the question of criticality on one hand, as some authors argue that information transmission is maximized near the critical values \cite{beggs, werner}, and the question of metastability in the other hand, as metastable dynamics seems to play an important role in the ability of the brain to process information \cite{deco, robert}. Another important interrogation is the topology of the networks, that is the question of the family of graph to which actual neural networks belong. While many complex real life networks, such as internet networks and social interacting species, have been extensively studied \cite{courgette}, leading to famous families of graphs such as Small World graphs, a consensus on the topology of neuronal networks has yet to be done. The model we consider here was originally introduced and studied in \cite{ferrari} on the infinite one-dimensional lattice. It was proven in this case that the model is subject to a phase-transition with respect to the leakage parameter $\gamma$, namely there exists a critical value $\gamma_c$ such that for any $\gamma < \gamma_c$ each neuron continues to spike infinitely many times with positive probability, while for $\gamma > \gamma_c$ the spiking activity stop with probability $1$. The same model was then studied in \cite{andre} on the finite one-dimensional lattice of arbitrary length and it was proven that for $\gamma$ small enough the renormalized time of extinction of the system is asymptotically  memory-less when the number of neurons grows, or in other words, asymptotically exponentially distributed. This feature is a crucial property of metastable dynamics as proposed in the seminal paper \cite{cassandro}, which was originally concerned with interacting particle physics. 

\vspace{0.4 cm}

Here we first prove that for the model defined on the finite one-dimensional lattice of arbitrary length this last result doesn't hold anymore if you take $\gamma$ big enough ($\gamma > 1$ is sufficient). Furthermore, for such a $\gamma$ the renormalized extinction time is asymptotically deterministic which makes it highly non-metastable. Then we turn ourself to the finite complete graph, and we prove that, on the contrary, the result of metastability holds for any positive $\gamma$ in this case. This introduces an interesting dichotomy with respect to metastability between loosely connected graphs, like the one-dimensional lattice, and highly connected ones, such as the complete graph.

\vspace{0.4 cm}

This paper is organized as follows. In Section \ref{def} we give a formal definition of the model, and we state the previously proven results as well as the two results that this article is concerned with. In Section \ref{reslattice} we give the proof of the result related to the one-dimensional lattice. In Section \ref{rescomplete} we give the proof of the result related to the complete graph.

\section{Definition and results}
\label{def}

The general framework of the model we consider is as follows. $V$ is a finite or countable set representing the neurons, and to each $i \in V$ is associated a subset $\V_i \subset V$ of \textit{presynaptic neurons}. For any $i \in V$ the set $\{ j \in V : i \in  \mathbb{V}_j\}$ is the set of the \textit{postsynaptic neurons} of neuron $i$. If you consider the elements of $V$ as vertices, and draw and edge from $i$ to $j$ whenever $i \in \V_j$, that is if you define $E = \{(i,j) : i \in \V_j\}$, then you obtain the graph of the network $G=(V,E)$. Each neuron $i$ has a \textit{membrane potential} evolving over time, represented by a stochastic process which takes its values in the set $\mathbb{N}$ of non-negative integers and which is denoted $(X_i(t))_{t \geq 0}$. Moreover each neuron is also associated with an homogeneous Poisson process on $\R_+$ with intensity $\gamma$, denoted $(N^{\dagger}_i(t))_{t \geq 0}$, which represents the \textit{leak times}. At any of these leak times the membrane potential of the neuron $i$ is reset to $0$. Finally an other point process on $\R_+$, denoted $(N_i(t))_{t \geq 0}$ and representing the \textit{spiking times}, is also associated to each neuron, whose rate at time $t$ is given by $\phi (X_i(t))$, $\phi$ being some function called the rate function. When a neuron spikes its membrane potential is reset to $0$ and the membrane potential of all of its postsynaptic neurons is increased by one. All the point processes involved are assumed to be mutually independent, and they determine entirely the evolution of the membrane potentials. See \cite{ferrari} for more details.

\vspace{0.4 cm}

In the present article we consider the special case in which the function $\phi$ is the indicator function defined for any $x \in \R$ by $\phi(x) = \mathbbm{1}_{x > 0}$. In this case the dynamic can be reformulated as follows. At any given time $t \geq 0$ a neuron $i \in V$ can either have its membrane potential strictly positive ($\mathbbm{1}_{X_i(t) > 0} = 1$) or not ($\mathbbm{1}_{X_i(t) > 0} = 0$), and depending on which is true we say that neuron $i$ is \textit{active} or \textit{quiescent} respectively. Whenever a neuron is active, it is associated with two independent exponential random clocks of parameter $1$ and $\gamma$, and it spikes or leaks respectively depending on which one rings first. Whenever a spike occurs, the neuron that spiked becomes quiescent, while its post-synaptic neurons become active if they weren't already. When a leak occurs, the neuron that leaks simply becomes quiescent.

\vspace{0.4 cm}

The distinction between active and quiescent neurons which arises from our choice for the function  $\phi$ is reminiscent of the distinction found in the classical Wilson-Cowan model (\cite{cowan}). It is a simplification which will prove itself to be  mathematically very convenient, and it has the advantage of putting the focus on the macroscopic phenomenons rather than on the details of the evolution of the membrane potentials.

\vspace{0.4 cm}

For any fixed $t \geq 0$ and $i \in V$ the quantity of interest is therefore given by $\mathbbm{1}_{X_i(t) > 0}$, as it corresponds to the infinitesimal spiking rate of neuron $i$ at time $t$. If for any $i \in V$ and $t \geq 0$ we write $\eta_i(t) = \mathbbm{1}_{X_i(t) > 0}$, and $\eta(t) = (\eta_j(t))_{j \in V}$, then the resulting process $(\eta(t))_{t \geq 0}$ is an interacting particle system, that is, a Markovian process taking value in $\{0,1\}^V$ (see \cite{ips}). Its infinitesimal generator (see \cite{ferrari}) is given by

\begin{equation} \label{generator} \L f (\eta) = \gamma \sum_{i \in V} \Big(f(\pi^{\dagger}_i(\eta)) - f(\eta)\Big) + \sum_{i \in V} \eta_i\Big(f(\pi_i(\eta)) - f(\eta)\Big),
\end{equation} where $f : \{0,1\}^{V} \mapsto \R$ is a cylinder function, $\gamma$ is a non-negative real number, and the $\pi^{\dagger}_i$'s and $\pi_i$'s are maps from $\{0,1\}^{V}$ to $\{0,1\}^{V}$ defined for any $i \in V$ as follows:

\[
    {\Big(\pi^{\dagger}_i(\eta)}\Big)_j= 
\begin{cases}
    0& \text{if } j = i,\\
    \eta_j & \text{otherwise},
\end{cases}
\] and

\[
    {\Big(\pi_i(\eta)}\Big)_j= 
\begin{cases}
    0 & \text{if } j = i,\\
    \max (\eta_i, \eta_j)  & \text{if } i \in \V_j,\\
    \eta_j & \text{otherwise}.
\end{cases}
\]

\vspace{0.4 cm}

The $\pi^{\dagger}_i$'s correspond to the leakage effect mentioned above, and the $\pi_i$'s correspond to the spikes.

\vspace{0.4 cm}

In the present article we consider two cases for the graph of interaction which are at the two ends of the spectrum with respect to connectivity: the case of a nearest-neighbors interaction on the one-dimensional lattice and the case in which the interaction is complete. As already noticed in \cite{ferrari}, the case of the nearest-neighbors interaction is a radical simplification, and it doesn't apply to the brain structure in general, but this kind of structure can be found in simpler nervous tissue, such as the retina (see \cite{braitenberg}). On the other hand, as actual biological neural network are usually connected—in the sense that there is always a path from a neuron to another—the complete interaction can be seen as some kind of mean-field approximation.

\vspace{0.4 cm}

The first case corresponds to the following assumption.

\begin{assumption} \label{latticeassumption}
$V = \Z$ and for all $i \in V$, $\V_i = \{i-1, i+1\}$.
\end{assumption}

It has been proven in \cite{ferrari} that this instantiation is subject to a phase transition. More precisely the following theorem was proven.

\begin{theorem} \label{thm:phasetransition}
Under Assumption \ref{latticeassumption} and assuming that $X_i(0) \geq 1$ for all $i \in V$, there exists a critical value $\gamma_c$ for the parameter $\gamma$, with $0 < \gamma_c < \infty$, such that for any $i \in \Z$

$$\P \Big( N_i([0,\infty[) \text{ } < \infty \Big) = 1 \text{ if } \gamma > \gamma_c$$

and

$$ \P \Big( N_i([0,\infty[) \text{ } = \infty \Big) > 0 \text{ if } \gamma < \gamma_c.$$

\end{theorem}

In other words there exists a critical value $\gamma_c$ for the leakage rate which is such that the process dies almost surely above it, and survive with positive probability below it.

\vspace{0.4 cm}
For this one-dimensional lattice instantiation let $\xi_i(t) \mydef \mathbbm{1}_{X_i(t) > 0}$ for any $i \in \Z$ and $t \geq 0$, so that $\big(\xi(t)\big)_{t \geq 0}$ is the interacting particle system taking value in $\{0,1\}^\Z$ whose dynamic is specified by (\ref{generator}). Now consider a finite version of the process $\big(\xi(t)\big)_{t \geq 0}$, where the neurons aren't anymore in the infinite lattice but only in a finite portion of it. Fix some integer $N \geq 0$ and let $V_N = \llbracket -N,N \rrbracket$ (where $\llbracket -N,N \rrbracket$ is a short-hand for $\Z \cap [-N,N]$). For any $i \in V_N$ let the set $\V_{N,i}$ be defined by

\[
    \V_{N,i}= 
\begin{cases}
    \{i-1,i+1\} & \text{if } i \in \llbracket -(N-1),N-1  \rrbracket,\\
    \{N-1\}  & \text{if } i = N,\\
    \{-(N-1)\}& \text{if } i = -N.
\end{cases}
\]

We write $\big(\xi_N(t)\big)_{t \geq 0}$ for the finite version of the lattice process, that is to say the process on $\{0,1\}^{2N+1}$ in which the neurons are indexed on $V_N$ and where the set of presynaptic neurons for neuron $i$ is given by $\V_{N,i}$. By convention, when no indication is given, the initial state of the process is the state where all neurons are active. Notice that in this finite setup the state in which all neurons are quiescent is an absorbing state, so that we can define the time of extinction of this finite process, which we denote $\tau_N$. Formally $$\tau_N \mydef \inf \{t \geq 0 : \xi_N(t)_i = 0 \text{ for all } i \in  \llbracket -N,N\rrbracket\}.$$

\vspace{0.4 cm}

It was proven in \cite{andre} that the following holds: 

\begin{theorem} \label{thm:metastabilityLattice}
There exists a $\gamma_c'$ satisfying $\gamma_c' \leq \gamma_c$ such that if $0 < \gamma < \gamma_c'$ then $\tau_N$ has finite expected value and

$$\frac{\tau_N}{\E (\tau_N)} \overset{\mathcal{D}}{\underset{N \rightarrow \infty}{\longrightarrow}} \mathcal{E} (1),$$ where $\E$ denotes the expectation, $\mathcal{D}$ denotes a convergence in distribution and $\mathcal{E} (1)$ an exponential random variable of mean $1$.
\end{theorem}

\vspace{0.4 cm}

In this article the first result obtained show that Theorem \ref{thm:metastabilityLattice} doesn't hold anymore if $\gamma$ is big enough. In fact the following theorem is proven.

\begin{theorem} \label{convone}
Suppose that $\gamma > 1$. Then the expectation of $\tau_N$ is finite and the following convergence holds

$$\frac{\tau_N}{\E (\tau_N)}  \overset{\P}{\underset{N \rightarrow \infty}{\longrightarrow}} 1,$$

where $\P$ denotes a convergence in probability.
\end{theorem}

\vspace{0.4 cm}

This result is in some sense symmetrical to Theorem \ref{thm:metastabilityLattice}. Indeed the later tells us that in a portion of the sub-critical region the time of extinction is asymptotically memory-less, which means that it is highly unpredictable: knowing that the process survived up to time $t$ doesn't give you any information about what should happen after time $t$. What we prove here is that in a  portion of the super-critical regime (indeed $\gamma_c < 1$ as proven in in \cite{ferrari}) the time of extinction is asymptotically constant, so that it is highly predictable.

\vspace{0.4 cm}

In the second part of this article we study the behaviour of the same stochastic process when the graph of the network is the complete graph of size $N$. 

\vspace{0.4 cm}

For some $N \in \N$ we let $V_N'= \llbracket 1, N\rrbracket$ and for any $i \in V_N'$ we let $V_{N,i}' = V_N' \setminus \{i\}$. Now write $\big(\zeta_N(t)\big)_{t \geq 0}$ for the equivalent of $\big(\xi_N(t)\big)_{t \geq 0}$ in the complete setting, that is the spiking rate process in which the neurons are indexed on $V_N'$ and where the set of presynaptic neurons for neuron $i$ is given by $V_{N,i}'$. In this instantiation every neuron is connected to one another so that when a single active neuron spikes, every other neuron becomes active. Let $\sigma_N$ denote the time of extinction of this finite process,  
$$\sigma_N \mydef \inf \{t \geq 0 : \zeta_N(t)_i = 0 \text{ for all } i \in  \llbracket 1,N\rrbracket\}.$$

\vspace{0.4 cm}

The second result we're aimed to prove is that for the complete graph instantiation the result of convergence of the renormalized time of extinction toward an exponential random variable that was proven in the lattice case for small $\gamma$ now holds for any positive $\gamma$. This is the object of the following theorem.

\vspace{0.4 cm}

\begin{theorem} \label{thm:convexpcompl}
For any $\gamma > 0$ the expectation of $\sigma_N$ is finite and the following convergence holds

$$\frac{\sigma_N}{\E (\sigma_N)} \overset{\mathcal{D}}{\underset{N \rightarrow \infty}{\longrightarrow}} \mathcal{E} (1).$$
\end{theorem}

\vspace{0.4 cm}

\section{Result on the one-dimensional lattice}
\label{reslattice}

In the following we will repeatedly identify the state space $\{0,1\}^\Z$ with $\mathcal{P}(\Z)$, the set of all subsets of $\Z$. Indeed any state $\eta$ of the process, belonging to $\{0,1\}^{\Z}$, can be seen as well as an element $A$ of $\mathcal{P}(\Z)$, writing $A = \{i \in \Z \text{ such that } \eta_i = 1\}$. For all the processes involved in the proof we will adopt the convention of writing the initial state as a superscript. For example we will write $(\xi^i(t))_{t \geq 0}$ and $(\xi^i_N(t))_{t \geq 0}$ to denote the infinite-lattice process and the finite version respectively, starting at time $0$ with only neuron $i$ active and all other neurons quiescent. Note that to avoid an overload in the notation we write $i$ instead of $\{i\}$, which is an abuse of notation. As stated in the previous section, in the absence of superscript the initial state is the state where all neurons are active, in particular $(\xi_N(t))_{t \geq 0}$ will denote the finite version of the process starting from the state $\llbracket -N,N\rrbracket$.

\vspace{0.4 cm}

The present section contains two propositions. The first one tells us that the time of extinction evolves asymptotically like a logarithm in $N$, and the second one tells us that the same is true for its expectation. Theorem \ref{convone} is an immediate consequence of these two propositions. Our proof use some of the ideas developed in \cite{liu} concerning the Harris contact process.

\begin{proposition} \label{convproba}
Suppose that $\gamma > 1$. Then there exists a constant $0<C<\infty$ depending on $\gamma$ such that the following convergence holds

$$ \frac{\tau_N}{\log (2N+1)} \overset{\P}{\underset{N \rightarrow \infty}{\longrightarrow}} C.$$
\end{proposition}

\begin{proof}

We define the following function

$$t \mapsto f(t) = \log \Big( \P \left( \xi^0 (t) \neq \emptyset \right)\Big).$$

We also define the following constant

$$ C' = -\sup_{s>0} \frac{f(s)}{s}.$$

\vspace{0.4 cm}

The first step is to show that the function $f$ is superadditive. For any $s,t \geq 0$ we have:

$$ \P\left(\xi^0 (t+s) \neq \emptyset \text{ } | \text{ } \xi^0(t) \neq \emptyset \right) \geq \P \left(\xi^0 (s) \neq \emptyset \right).$$

\vspace{0.4 cm}

This last inequality follows from Markov property and the fact that having a higher number of active neurons in the initial configuration implies having a higher probability of being alive for the process at any given time $s$ (see Proposition 4.2 in \cite{andre}). Moreover if $\xi^0(t) \neq \emptyset$ then $|\xi^0(t)| \geq 1 = |\xi^0(0)|$. Furthermore it can be rewritten as follows

$$ \P\left(\xi^0 (t+s) \neq \emptyset \right) \geq \P \left(\xi^0 (t) \neq \emptyset \right) \P \left(\xi^0 (s) \neq \emptyset \right),$$ 
and taking the log gives the superaddtivity we are looking for. Now from a well-known result about superadditive functions, sometimes called the Fekete lemma \cite{fekete}, we get the following convergence

\begin{equation}
\label{superaddconv}
\frac{f(t)}{t} {\underset{t \rightarrow \infty}{\longrightarrow}} - C'.
\end{equation}

\vspace{0.4 cm}

For any $t > 0$ we also have

\begin{equation}
\label{expbound}
\P \left( \xi^0 (t) \neq \emptyset \right) \leq e^{-C't}.
\end{equation}

\vspace{0.4 cm}

Notice that while it is clear that $0 \leq C' < \infty$, it is not obvious that $C' > 0$. We show that it is the case using a coupling with the following  branching process. At time $0$ there is a single individual. Two independent exponential random clocks of parameter $1$ and $\gamma$ respectively are attached to this individual. If the rate $\gamma$ clock rings before the other one, then the individual dies. In the contrary case the individual is replaced by two other individuals. Every new individual gets his two own independent exponential clocks of parameter $1$ and $\gamma$ and so on. We denote by $(Z_t)_{t \geq 0}$ the process corresponding to the number of individuals of the population along the time. Note that by hypothesis we have $Z_0 = 1$. The expectation at time $t$ can be explicitly computed (see for example chapter 8 in \cite{schinazi}). For any $\gamma \geq 0$ and $t \geq 0$ we have

\begin{equation}
\label{expectation}
 \E \left(Z_t \right) = e^{- (\gamma - 1)t}.
\end{equation}

\vspace{0.4 cm}

The coupling is done as follows, at time $0$ the only active neuron in $\xi^0(0)$ is coupled with the only individual in $Z(0)$. By this we mean that if this neuron becomes quiescent then the individual dies, and if the neuron spikes, then the individual is replaced by two new individuals. When a spike occurs, there is three possibilities: two neurons are activated, one neuron is activated, or no neuron is activated, depending on how much neighbours are quiescent. In any of these cases each activated neuron is coupled with a newborn individuals, and any supernumerary newborn individual is given his own independent exponential clocks. At any time $t \geq 0$ we obviously have $|\xi^0_t| \leq Z_t$. Using (\ref{expectation}) and Markov inequality it follows that

$$\P \left( \xi^0(t) \neq \emptyset \right) \leq \P \left( Z_t \geq 1 \right) \leq e^{-(\gamma - 1)t}.$$

\vspace{0.4 cm}

Then we take the log and divide by $t$ in the previous inequality and we obtain at the limit that $C' \geq \gamma - 1$, and from the assumption that $\gamma > 1$ we get $C' > 0$.

\vspace{0.4 cm}

Let us break the suspense and already reveal that the constant $C$ we are looking for is actually simply the inverse of $C'$. Therefore in order to prove our result we are going to prove that for any $\epsilon > 0$ we have the two following convergences

\begin{equation}
\label{pluseps}
    \P \left( \frac{\tau_N}{\log(2N+1)} - \frac{1}{C'} > \epsilon \right) {\underset{N \rightarrow \infty}{\longrightarrow}} 0,
\end{equation}
and

\begin{equation}
\label{moinseps}
    \P \left( \frac{\tau_N}{\log(2N+1)} - \frac{1}{C'} < -\epsilon \right) {\underset{N \rightarrow \infty}{\longrightarrow}} 0.
\end{equation}

\vspace{0.4 cm}

Let us start with (\ref{pluseps}), which is the easiest part. We remark that our process is \textit{additive} in the sense that $\xi_N(t) = \cup_{i \in \llbracket -N,N \rrbracket} \xi^i_N(t)$. See \cite{ferrari} for details, it is also an immediate consequence of the graphical construction proposed in \cite{andre}. Using additivity and inequality (\ref{expbound}) we get

\begin{equation} \label{expboundN}
\P \left( \xi_N(t) \neq \emptyset \right) \leq (2N + 1) \P \left( \xi^0(t) \neq \emptyset \right) \leq (2N+1) e^{-C't}.
\end{equation}

\vspace{0.4 cm}

Now, for any $\epsilon > 0$, if you let $t = (\frac{1}{C'} + \epsilon)\log(2N+1)$ then the following holds

$$\P \left( \frac{\tau_N}{\log(2N+1)} - \frac{1}{C'} > \epsilon \right) = P \left( \xi_N(t) \neq \emptyset \right) \leq e^{-C'\epsilon \log(2N + 1)}.$$

Then the fact that $C' > 0$ ensures us that the term on the right-hand side of the inequality goes to $0$ as $N$ diverges, which proves (\ref{pluseps}).

\vspace{0.4 cm}

It remains to prove (\ref{moinseps}). If for some $N \in \N^*$ we take $t = \left(\frac{1}{C'} - \epsilon \right) \log (2N + 1)$, then we can write

$$ \P \left( \frac{\tau_N}{\log(2N+1)} - \frac{1}{C'} < -\epsilon \right) = \P \left( \xi_N(t) = \emptyset \right),$$
so that it suffices to show that the right-hand side converges to $0$ for this choice of $t$ as $N$ goes to $\infty$. For reasons that will become clear in a moment we will actually write $t = \frac{1}{C'} \left(1 - \epsilon' \right) \log (2N + 1)$, with $\epsilon' = C' \epsilon$.

\vspace{0.4 cm}

From (\ref{superaddconv}) (and from the fact that $C'> 0$) we get that for any $\epsilon > 0$ and for big enough $t$

$$ \frac{f(t)}{t} \geq - (1+\epsilon') C',$$
which can be written

$$\P \left( \xi^0_t = \emptyset \right) \leq 1 - e^{- (1+\epsilon') C't}.$$

Therefore, with $t = \frac{1}{C'} \left(1 - \epsilon' \right) \log (2N + 1)$ and $N$ big enough we have

\begin{equation} \label{boundxi0}
\P \left( \xi^0_t = \emptyset \right) \leq 1 - \frac{1}{(2N+1)^{1-\epsilon'^2}}.
\end{equation}

\vspace{0.4 cm}

Now for any $k \in \Z$ we define

$$F_k \text{ } \mydef \text{ }  \llbracket(2k-1)K\log(2N+1), (2k+1)K\log(2N+1)\rrbracket,$$
where $K$ is some constant depending on $N$ whose value will be chosen later in order for $K\log(2N+1)$ to be an integer. We then consider a modification of the process $(\xi_N(t))_{t \geq 0}$ where all neurons at the border of one of the sub-windows $F_k$ defined above (i.e. all neurons indexed by $(2k+1)K\log(2N+1)$ for some $k \in \Z$) are fixed in quiescent state and therefore are never allowed to spike. This modified process is denoted $(\widetilde{\xi}_N(t))_{t \geq 0}$. We also define the following configuration

$$ A_N \text{ } \mydef \text{ } \left\{ 2kK\log(2N+1) \text{ for } k \in \Z \cap \left[-\frac{N}{2K\log(2N+1)}, \frac{N}{2K\log(2N+1)}\right] \right\}.$$

\vspace{0.4 cm}

Notice that the fact that the neurons at the borders of the windows $F_k$ are never allowed to spike makes the evolution of $(\widetilde{\xi}_N(t))_{t \geq 0}$ independent from one window to another. Moreover notice that the integers belonging to $A_N$ are all at the center of one of these windows. 

\vspace{0.4 cm}

Now for any $t \geq 0$ we define $r_t \mydef \max \xi^0_t$. Considering the spiking process $(\xi_t)_{t \geq 0}$ with no leaking it is easy to see that the right edge $r_t$ can be coupled with an homogeneous Poisson process of parameter $1$, that we denote $(M(t))_{t \geq 0}$, in such a way that for any $m \geq 0$

$$ \P \left( \sup_{s \leq t} r_s \geq m \right) \leq \P \Big( M(t) \geq m \Big).$$

\vspace{0.4 cm}

We have

$$ \E \left( e^{M(t)} \right) = e^{t(e - 1)},$$

so taking the exponential, using Markov inequality and taking $m=K't$ (where $K'$ is some constant that we are going to fix in a moment) we get

\begin{align*}
 \P \left( \sup_{s \leq t} r_s \geq K't\right) & \leq e^{t(e - 1 - K')}\\
 & \leq e^{t(2 - K')},
\end{align*}

where in the last inequality we simply used the fact that $e - 1 < 2$.

\vspace{0.4 cm}

Now taking again $t = \frac{1}{C'}\left(1 - \epsilon' \right) \log (2N + 1)$ and $K' = 2(1 + C')$  we get

$$ \P \left( \sup_{s \leq t} r_s \geq m \right) \leq e^{-2(1-\epsilon') \log(2N+1)},$$

and assuming without loss of generality that $\epsilon' < \frac{1}{2}$ we get

\begin{equation}
\label{boundEt}
\P \left( \sup_{s \leq t} r_s \geq m \right) \leq \frac{1}{2N+1}.
\end{equation}

\vspace{0.4 cm}

It is now possible to fix the value of the constant $K$ we introduced earlier. We take

$$ K = \inf \left\{ x \in \R \text{ such that } x \geq \frac{K'}{C'} \text{ and } x \log(2N+1) \in \N \right\}.$$

In other words we take $K$ equal to $\frac{K'}{C'}$ and then enlarge it slightly in order for $K\log(2N+1)$ to be an integer. We also define the following event

$$ E_t \text{ } \mydef \text{ } \Big\{ \xi^0_s \text{ doesn't escape from } \llbracket-K\log(2N+1), \ldots, K\log(2N+1)\rrbracket \text{ for any } s \leq t\Big\}.$$

Now taking $N$ large enough and $t = \frac{1}{C'}\left(1 - \epsilon' \right) \log (2N + 1)$ we have

\begin{align*}
\P \left( \xi_N(t) = \emptyset \right) &\leq \P \left( \widetilde{\xi}^{A_N}_N (t) = \emptyset \right)\\
& = \P \left( \widetilde{\xi}^{0}_N(t) = \emptyset \right)^{(2N+1) / (2K\log(2N+1))}\\
& \leq \Big(  \P \left( \widetilde{\xi}^0_N (t) = \emptyset \cap E_t\right) + \P \left( E_t^c \right) \Big)^{(2N+1) / (2K\log(2N+1))}\\
& \leq \Big(  \P \left( \xi^0 (t) = \emptyset\right) + \P \left( E_t^c \right) \Big)^{(2N+1) / (2K\log(2N+1))}\\
& \leq \left( 1 - \left(\frac{1}{(2N+1)^{1 - \epsilon'^2}} - \frac{2}{2N+1} \right) \right)^{(2N+1) / (2K\log(2N+1))}.\\
\end{align*}

\vspace{0.4 cm}

To obtain the inequality above we used (\ref{boundxi0}) and the fact that inside $E_t$ the process $(\widetilde{\xi}^0_N(s))_{0 \leq s \leq t}$ evolves just like $(\xi^0_s)_{0 \leq s \leq t}$, which allows us to bound $\P \left( \xi^0_N (t) = \emptyset\right)$, and we used (\ref{boundEt}) to bound $\P \left( E_t^c \right)$.

\vspace{0.4 cm}

Finally we let 

$$ a_N = \frac{1}{(2N+1)^{1 - \epsilon'^2}} - \frac{2}{2N+1},$$

and 

$$ b_N = \frac{2N+1}{2K\log(2N+1)},$$

so that the last bound can be written $(1 - a_N)^{b_N}$. Then $$(1 - a_N)^{b_N} = e^{b_N \log (1 - a_N)} \leq e^{-b_N a_N},$$ and since $a_Nb_N \underset{N \rightarrow \infty}{\longrightarrow} \infty$, it follows that $e^{-b_N a_N}$ goes to $0$ as $N$ goes to $\infty$.

\end{proof}

\vspace{0.4 cm}

The last step consists in showing that the same convergence holds for the expectation, which is the object of the following proposition.

\begin{proposition} \label{convexpect}
Suppose that $\gamma > 1$. Then the expectation of $\tau_N$ is finite and the following convergence holds

$$ \frac{\E \left(\tau_N\right)}{\log (2N+1)} \underset{N \rightarrow \infty}{\longrightarrow} C,$$ where $C$ is the same constant as in Proposition \ref{convproba}.
\end{proposition}

\begin{proof}
It is well-known that the fact that a sequence of random variables $(X_n)_{n \in \N}$ converges in probability to some random variable $X$ doesn't necessarily implies that $\E(X_n) \underset{n \rightarrow \infty}{\longrightarrow} \E(X)$. Nonetheless this implication holds true with the additional assumption that the sequence is uniformly integrable (see for example Theorem 5.5.2 in \cite{tande} page 259), i.e. if the following holds

$$ \lim_{M \rightarrow \infty} \left( \sup_{n \in \N} \E \Big( |X_n| \mathbbm{1}_{\{|X_n| > M\}} \Big)\right) = 0.$$

\vspace{0.4 cm}

It is therefore sufficient to show that $\big(\tau_N/\log(2N+1)\big)_{N \in \N^*}$ is uniformly integrable, and the result will follows from Proposition \ref{convproba}. For some $M > 0$ and some $N \in \N^*$ it is easy to see that we have the following

$$ \E \left( \frac{\tau_N}{\log (2N+1)} \mathbbm{1}_{\{\frac{\tau_N}{\log (2N+1)} > M\}} \right) = \int_0^\infty \P \left( \frac{\tau_N}{\log (2N+1)} > \max(t,M) \right) dt.$$ 

\vspace{0.4 cm}
Now using inequality (\ref{expboundN}) and the previously proven fact that $C'>0$ when $\gamma > 1$ we have the following

\begin{align*}
&\int_0^\infty \P \left( \frac{\tau_N}{\log (2N+1)} > \max(t,M) \right) dt\\
&= \int_0^M \P \left( \frac{\tau_N}{\log (2N+1)} > M \right) dt + \int_M^\infty \P \left( \frac{\tau_N}{\log (2N+1)} > t \right) dt\\
&\leq (2N +1) \left[ \int_0^M e^{-C'\log (2N+1)M} dt + \int_M^\infty e^{-C'\log (2N+1)t} dt \right]\\
&= (2N +1)^{1 - C'M} \left[ M + \frac{1}{C'\log(2N+1)}\right],
\end{align*}where $C'$ is the same constant as in the previous proof. Without loss of generality we assume that $M > \frac{1}{C'}$, so that the bound above is decreasing in $N$, from what we get

\begin{equation} \label{supbound}
\sup_{n \in \N^*} \E \left( \frac{\tau_N}{\log (2N+1)} \mathbbm{1}_{\{\frac{\tau_N}{\log (2N+1)} > M\}}\right) \leq 3^{1 - C'M} \left[ M + \frac{1}{C'\log(3)}\right].
\end{equation}

\vspace{0.4 cm}

Finally the right-hand side of inequality (\ref{supbound}) goes to $0$ when $M$ goes to $\infty$, so that the uniform integrability is proven. Notice that moreover the finiteness of the expectation is an immediate consequence of uniform integrability.
\end{proof}

\vspace{0.4 cm}

\section{The complete graph case}
\label{rescomplete}

Recall that $(\zeta_N (t))_{t \geq 0}$ denotes the analog of $(\xi_N(t))_{t \geq 0}$ in the complete graph setting, and that $\sigma_N$ denotes its time of extinction. We adopt the same conventions as in the previous section regarding the notation, identifying $\{0,1\}^{\llbracket 1,N\rrbracket}$ with $\mathcal{P}(\llbracket 1,N\rrbracket)$ and indicating the initial state as a superscript when it is different from the whole space $\llbracket 1,N\rrbracket$. This convention applies to the time of extinction as well. Notice that in the complete graph setting it doesn't really make sense to consider an infinite process anymore, as any neuron would be activated by a spike infinitely many times in any time interval for example, moreover the minimal conditions for the existence of an infinitesimal generator are not satisfied. Luckily we wont need any infinite process in the course of the proof.

\vspace{0.4 cm}

The Markov process $(\zeta_N (t))_{t \geq 0}$ has three communicating classes, which are visited in increasing order. The first one contains only the "all one" state, which is left definitively after the first spike/leak, and the third one contains only the "all zero" state, which is an absorbing state. The second class contains all the other states. As the absorbing state is accessible from the second class, standard results on irreducibility tell us that the time of extinction $\sigma_N$ is almost surely finite. It implies that $\P (\sigma_N > t)$ converges to $0$ as $t$ diverges.

\vspace{0.4 cm}

Moreover it is clear from the definition of $(\zeta_N (t))_{t \geq 0}$ that $\P (\sigma_N > t)$ is a continuous and decreasing function of $t$. Therefore we can define $\beta_N$ to be the unique value in $\R^+$ such that

$$ \P (\sigma_N > \beta_N) = e^{-1}.$$

\vspace{0.4 cm}

The main ingredient of the proof is the Proposition \ref{prop:difftozero} below. Theorem \ref{thm:convexpcompl} then follows from Corollary \ref{cor:convbeta} and Proposition \ref{prop:betaexpequi}, which tell us respectively that $\sigma_N / \beta_N$ converges in distribution to an exponential random variable of mean $1$ and that $\E\left( \sigma_N \right) \underset{N \infty}{\sim} \beta_N$.

\vspace{0.4 cm}

\begin{proposition} \label{prop:difftozero}
For any $\gamma > 0$ and for any $s,t \geq 0$ the following holds

\begin{align} \label{result}
\lim_{N \rightarrow \infty} \left| \P\left( \frac{\sigma_N}{\beta_N} > s + t \right) - \P \left( \frac{\sigma_N}{\beta_N} > s \right)\P\left( \frac{\sigma_N}{\beta_N} > t \right) \right| = 0.
\end{align}
\end{proposition}

\vspace{0.4 cm}
\begin{proof}

In our complete graph setting there is no spatial dependence between the neurons like in the lattice setting, so that the law of the time of extinction is impacted by the initial state only through its cardinal. In particular, for any $1 \leq k \leq N$ and any $A \in \mathcal{P}(\llbracket 1,N \rrbracket)$ of size $k$ , $\sigma^A_N$ has the same law as $\sigma_N^{\llbracket 1, k\rrbracket}$. Using the Markov property and this last remark we have

\begin{align*}
&\left| \P\left( \frac{\sigma_N}{\beta_N} > s + t \right) - \P \left( \frac{\sigma_N}{\beta_N} > s \right)\P\left( \frac{\sigma_N}{\beta_N} > t \right) \right| \\
&= \sum_{k=1}^{N} \left|\mathbb{P}\left(\frac{\sigma_N^{\llbracket 1,k\rrbracket}}{\beta_N}>t\right)-\mathbb{P}\left(\frac{\sigma_N}{\beta_N}>t\right)\right| \cdot \mathbb{P}\big(|\zeta_N(\beta_N s)| = k \big) \\
&\leq \sum_{k=1}^{\left\lfloor \frac{N}{2} \right\rfloor} \mathbb{P}\left(|\zeta_N( \beta_N s)| = k \right) + \sum_{k=\left\lceil \frac{N}{2} \right\rceil}^{N} \left|\mathbb{P}\left(\frac{\sigma_N^{\llbracket 1,k\rrbracket}}{\beta_N}>t\right)-\mathbb{P}\left(\frac{\sigma_N}{\beta_N}>t\right)\right|
\end{align*}

\vspace{0.4 cm}

We fix $\epsilon > 0$. In order to prove the desired result we show  that we can find $N$ big enough such that

\begin{equation}
\label{sumbign}
\sum_{k=\left\lceil \frac{N}{2} \right\rceil}^{N} \left|\mathbb{P}\left(\frac{\sigma_N^{\llbracket 1,k\rrbracket}}{\beta_N}>t\right)-\mathbb{P}\left(\frac{\sigma_N}{\beta_N}>t\right)\right| < \epsilon,
\end{equation}
and

\begin{equation}
\label{sumsmalln}
\sum_{k=1}^{\left\lfloor \frac{N}{2} \right\rfloor} \mathbb{P}\left(|\zeta_N(\beta_N s)| = k \right) < \epsilon.
\end{equation}

\vspace{0.4 cm}

We start with (\ref{sumbign}), which is the easiest part. For any $\left\lceil \frac{N}{2} \right\rceil \leq k \leq N$ we denote by $E_k$ the event in which every active neuron in the process starting from $\llbracket 1, k\rrbracket$ becomes quiescent before any of them is affected by a spike. Notice that on the complementary event $E^c_k$ there is a spike affecting the process at some point, and that the process starting from  $\llbracket1, k\rrbracket$ and the process starting from $\llbracket1, N\rrbracket$ become equal at this point (both will be in the state $\llbracket 1,N \rrbracket \setminus \{i\}$, where $i$ is the neuron that just spiked). It follows that we have 

$$\mathbb{P} \left(\frac{\sigma_N^{\llbracket 1,k\rrbracket}}{\beta_N}>t \ \Big| \ E^c_k \right) = \mathbb{P} \left(\frac{\sigma_N}{\beta_N} > t \ \Big| \ E^c_k\right).$$

\vspace{0.4 cm}

From this we get that, for any $\left\lceil \frac{N}{2} \right\rceil \leq k \leq N$, 

\begin{align*}
&\left|\mathbb{P}\left(\frac{\sigma_N^{\llbracket 1,k\rrbracket}}{\beta_N}>t\right)-\mathbb{P}\left(\frac{\sigma_N}{\beta_N}>t\right)\right| \\
&= \left|\mathbb{P} \left(\frac{\sigma_N^{\llbracket 1,k\rrbracket}}{\beta_N}>t \ \Big| \ E_k \right) - \mathbb{P} \left(\frac{\sigma_N}{\beta_N} > t \ \Big| \ E_k\right)\right| \cdot \mathbb{P}(E_k) \\
&\leq \mathbb{P}(E_k).
\end{align*}

\vspace{0.4 cm}

Now since $\mathbb{P}(E_k) = \left( \frac{\gamma}{1 + \gamma}\right)^k \leq \left(\frac{\gamma}{1 + \gamma}\right)^{\frac{N}{2}}$ for any $k \geq \left\lceil \frac{N}{2} \right\rceil$ we get the following bound

\begin{align} \label{eq:bound1}
\sum_{k=\left\lceil \frac{N}{2} \right\rceil}^{N} \left|\mathbb{P}\left(\frac{\sigma_N^{\llbracket 1,k\rrbracket}}{\beta_N}>t\right)-\mathbb{P}\left(\frac{\sigma_N}{\beta_N}>t\right)\right| \leq \frac{N}{2} \left(\frac{\gamma}{1 + \gamma}\right)^{\frac{N}{2}},
\end{align}

\vspace{0.4 cm}
which goes to zero as $N$ goes to infinity, so that we can find some $N_1$ such that (\ref{sumbign}) is satisfied for any $N \geq N_1$.

\vspace{0.4 cm}

We now turn to (\ref{sumsmalln}). We will use a coupling. Let $(\widetilde{\zeta}_N(t))_{t \geq 0}$ be the process defined as follows. For $t < \sigma_N$ the process is simply equal to $\zeta_N(t)$. Now suppose $i \in \llbracket 1,N \rrbracket$ is the last neuron active in $(\zeta_N(t))_{t \geq 0}$ before the extinction at $\sigma_N$, then when $i$ leaks in $(\zeta_N(t))_{t \geq 0}$, instead of leaking in $(\widetilde{\zeta}_N(t))_{t \geq 0}$ it spikes. Then the dynamic of $(\widetilde{\zeta}_N(t))_{t \geq 0}$ is the same as the dynamic of the process $(\zeta_N(t))_{t \geq 0}$ with the only difference that whenever there is only one neuron remaining it doesn't leak, and spike at rate $1 + \gamma$, so that there is no extinction for this  stochastic process. The reason for asking that the last neuron spike at rate $1 + \gamma$ instead of simply $1$ is that if $X$ and $Y$ are two independent random variable exponentially distributed with rate $1$ and $\gamma$ respectively, then one can easily compute that the law of $X$ conditioned on $X \leq Y$ is an exponential distribution of rate $1 + \gamma$. Now for any $1 \leq k \leq \left\lfloor \frac{N}{2} \right\rfloor$, we have

\begin{equation}\label{eq:inf}
\begin{split}
\P \Big( |\zeta_N(\beta_N s)| = k \Big) &= \P \Big( |\zeta_N (\beta_N s)| = k \ \big| \ \sigma_N > \beta_N s \Big) \P \Big( \sigma_N > \beta_N s \Big)\\
&= \P \Big( |\widetilde{\zeta}_N (\beta_N s)| = k \ \big| \ \sigma_N > \beta_N s \Big) \P \Big( \sigma_N > \beta_N s \Big)\\
&\leq \P \Big( |\widetilde{\zeta}_N(\beta_N s)| = k \Big),
\end{split}
\end{equation}

so that it will be sufficient to prove that for $N$ big enough we have

\begin{equation}
\label{sumsmallntilde}
\sum_{k=1}^{\left\lfloor \frac{N}{2} \right\rfloor} \mathbb{P}\left(|\widetilde{\zeta}_N (\beta_N s)| = k \right) < \epsilon.
\end{equation}

\vspace{0.4 cm}

One can easily see that the process $(|\widetilde{\zeta}_N(t)|)_{t \geq 0}$, that is to say the process counting the number of particles at any time $t$, is a Markov jump process taking value in $\llbracket 1,N\rrbracket$, which transition diagram is given in Figure \ref{fig:markovchain}.

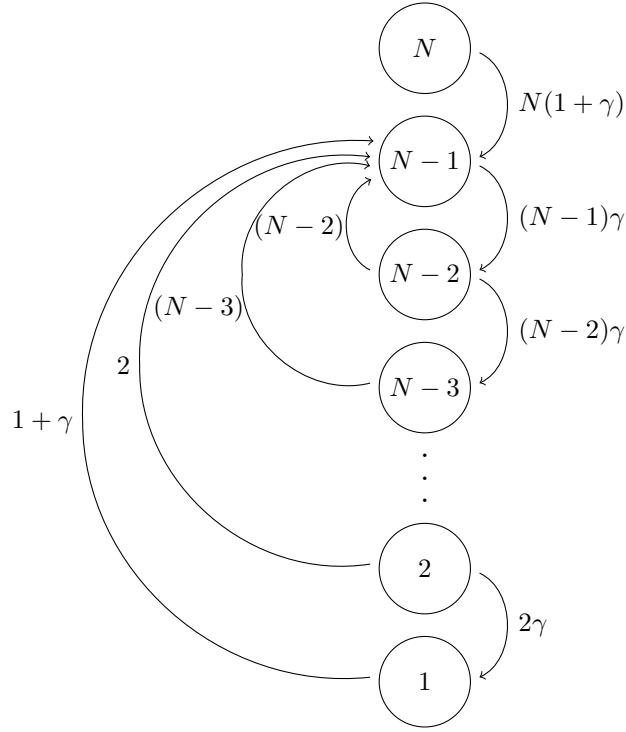
\begin{figure}[ht]
\begin{center}

\begin{tikzpicture}[scale=0.3]

\draw (0,15) circle (2) node[] {$N$}; 
\draw (0,10) circle (2) node[] {$N-1$}; 
\draw (0,5) circle (2) node[] {$N-2$};
\draw (0,0) circle (2) node[] {$N-3$};
\draw (0,-8) circle (2) node[] {$2$};
\draw (0,-13) circle (2) node[] {$1$};

\draw[->] (2.4,14.8) to [bend left=65] (2.4,10.2) ;
\draw[->] (2.4,9.8) to [bend left=65] (2.4,5.2) ;
\draw[->] (2.4,4.8) to [bend left=65] (2.4,0.2);
\draw[->] (2.4,-8.2) to [bend left=65] (2.4,-12.8);

\draw[->] (-2.4,5.2) to [bend left=65] (-2.35,9.2) ;
\draw[->] (-2.4,0.2) [bend left=55] to (-8,5) to [bend left=55] (-2.4,9.8) ;
\draw[->] (-2.4,-7.8) [bend left=50] to (-12.5,1) to [bend left=50] (-2.4,10.2) ;
\draw[->] (-2.4,-12.8) [bend left=48] to (-15,-1.5) to [bend left=48] (-2.25,10.9) ;

\draw[] (3.7,12.5) node[right] {$N(1+\gamma)$};
\draw[] (3.7,7.5) node[right] {$(N-1)\gamma$};
\draw[] (3.7,2.5) node[right] {$(N-2)\gamma$};
\draw[] (3.7,-10.5) node[right] {$2 \gamma$};

\draw[] (-3.1,7.1) node[left] {$(N-2)$};
\draw[] (-7.5,3.5) node[left] {$(N-3)$};
\draw[] (-12.5,1) node[left] {$2$};
\draw[] (-15,-1.5) node[left] {$1+\gamma$};

\draw[] (0,-3) node[] {\Large{.}};
\draw[] (0,-4) node[] {\Large{.}};
\draw[] (0,-5) node[] {\Large{.}};

\end{tikzpicture}
 \end{center}

\caption{The transition diagram of the Markov jump process $(\widetilde{\zeta}_N(t))_{t \geq 0}$. \label{fig:markovchain}} 
\end{figure}

\vspace{0.4 cm}

We would like to compute an invariant measure for this chain, so that we need to solve the following equation for $\mu$.

$$\mu Q = 0,$$

where $Q$ is the transition intensities matrix, given by

\vspace{0.4 cm}

\resizebox{0.96\textwidth}{!}{
$Q=\begin{pmatrix}
-(1+\gamma) & 0 & 0 & \hspace{1cm} & \cdots & \hspace{1cm} & 0 & 0 & 1+\gamma &0\\
2 \gamma & -2(1+\gamma) & 0 &  & \cdots & &0 & 0 & 2&0\\
0 & 3 \gamma  & -3(1+ \gamma) & & \cdots & & 0 & 0 & 3&0\\
\vdots & \vdots & \vdots & & \cdots & & \vdots & \vdots & \vdots & \vdots\\
0 & 0 & 0 & & \cdots & & (N-2)\gamma & -(N-2)(1+\gamma) & N-2 &0\\
0 & 0 & 0 & & \cdots & & 0 & (N-1)\gamma & -(N-1)\gamma &0\\
0 & 0 & 0 & & \cdots & & 0 & 0 & N(1+\gamma) & -N(1+\gamma) \\
\end{pmatrix}.$}

\vspace{0.4 cm}

This is the same as the following system of linear equations.

\begin{equation}
\begin{cases}
\label{linearsystem}
-(1+\gamma)\mu_1 + 2\gamma \mu_2 =0,\\
-2(1+\gamma )\mu_2 + 3\gamma \mu_3 = 0,\\
\ \ \ \ \ \ \ \ \ \ \ \ \ \ \vdots\\
-(N-2)(1+ \gamma) \mu_{N-2} + (N-1)\gamma \mu_{N-1} =0,\\
(1+\gamma)\mu_1 +2\mu_2 +3\mu_3 +  \ldots +(N-2)\mu_{N-2}-(N-1)\gamma\mu_{N-1} + N(1+\gamma)\mu_N  =0, \\
-N(1+\gamma)\mu_N=0.
\end{cases} 
\end{equation}

\vspace{0.4 cm}

Solving the system from top to bottom, we get that, for any $n \in \llbracket 2,N-1\rrbracket$ 

\begin{equation} \label{eq:mun}
\mu_n  = \frac{(1+\gamma)^{n-1}}{n \gamma^{n-1}}\mu_1.
\end{equation}

\vspace{0.4 cm}

Moreover $\mu_N = 0$. Now from (\ref{eq:mun}) and from the fact that the elements of $\mu$ need to sum up to $1$ we get

$$\mu_1 = \left(\sum_{n=1}^{N-1} \frac{(1+\gamma)^{n-1}}{n \gamma^{n-1}} \right)^{-1}.$$

\vspace{0.4 cm}

Disregarding all terms in the sum but the last one we obtain the following bound

$$\mu_1 < N\left(\frac{1+\gamma}{\gamma}\right)^{2-N}.$$

Hence, for any $1 \leq k \leq \left\lfloor \frac{N}{2} \right\rfloor$

$$\mu_k < N\left(\frac{1+\gamma}{\gamma}\right)^{k+1-N}.$$

\vspace{0.4 cm}

As a consequence,

\begin{equation} \label{eq:boundsuminv}
\sum_{k=1}^{\left\lfloor \frac{N}{2} \right\rfloor} \mu_k < \frac{N^2}{2}\left(\frac{1+ \gamma}{\gamma}\right)^{1-\frac{N}{2}}.
\end{equation}

\vspace{0.4 cm}

Let $\big(\widetilde{\zeta}_N^{\mu}(t)\big)_{t \geq 0}$ denote the process whose initial state is chosen according to the invariant measure $\mu$. By this we mean that a value $k \in \llbracket 1,N \rrbracket$ is sorted according to the invariant measure, and that the process then start from the initial state $\llbracket 1,k \rrbracket$. For any $t \geq 0$ we have the following inequality:

\begin{equation} \label{eq:suminv}
\mu_{N-1} \sum_{k=1}^{\left\lfloor \frac{N}{2} \right\rfloor} \mathbb{P}\left(\left|\widetilde{\zeta}^{\llbracket 1,N-1\rrbracket}_N(t)\right| = k \right) \leq \sum_{k=1}^{\left\lfloor \frac{N}{2} \right\rfloor} \mathbb{P}\left(|\widetilde{\zeta}^{\mu}_N(t)| = k \right).
\end{equation}

\vspace{0.4 cm}

This last inequality gets us closer to our goal but it is still not exactly what we need as the left hand side involves the process starting from $\llbracket 1,N-1\rrbracket$ while we would like it to start from the full initial configuration $\llbracket 1,N\rrbracket$. This little problem is solved as follows. Let $T_N$ be the time of the first jump of the process $(\widetilde{\zeta}_N(t))_{t \geq 0}$, that is to say

$$ T_N \mydef \inf \{ t \geq 0 : |\widetilde{\zeta}_N(t)| \neq N\}.$$

\vspace{0.4 cm}

Then, for any $t \geq 0$ and $1 \leq k \leq \left\lfloor \frac{N}{2} \right\rfloor$, the following holds

\begin{equation} \label{NminusonetoN}
\P \left( \left|\widetilde{\zeta}_N(t)\right| = k\right) = \P \left( \left|\widetilde{\zeta}^{\llbracket 1,N-1\rrbracket}_N\left((t - T_N)^+\right)\right| = k\right),
\end{equation}

where $(t - T_N)^+$ stands for $\max(0, t - T_N)$. This last inequality is obtained from Markov property and the fact that whenever $T_N > t$ the events we look at are both of probability $0$ for the $k$ we consider (assuming that $N \geq 3$).

\vspace{0.4 cm}

Now from (\ref{eq:boundsuminv}), (\ref{eq:suminv})  and (\ref{NminusonetoN}), we obtain

\begin{align*}
\sum_{k=1}^{\left\lfloor \frac{N}{2} \right\rfloor} \mathbb{P}\left(\left|\widetilde{\zeta}_N (\beta_N s)\right| = k \right) &= \sum_{k=1}^{\left\lfloor \frac{N}{2} \right\rfloor} \mathbb{P}\left(\left|\widetilde{\zeta}^{\llbracket 1,N-1\rrbracket}_N\left((\beta_N s - T_N )^+\right)\right| = k \right)\\
&\leq \frac{1}{\mu_{N-1}} \sum_{k=1}^{\left\lfloor \frac{N}{2} \right\rfloor} \mathbb{P}\left(\left|\widetilde{\zeta}^{\mu}_N\left((\beta_N s - T_N )^+\right)\right| = k \right)\\
&\leq \frac{1}{\mu_{N-1}} \frac{N^2}{2}\left(\frac{1+ \gamma}{\gamma}\right)^{1-\frac{N}{2}}.
\end{align*}

\vspace{0.4 cm}

Moreover from the penultimate line of (\ref{linearsystem}) we get $$ \mu_{N-1}  > \frac{ \mu_1 +\mu_2 +\mu_3 + \ldots +\mu_{N-2}}{\gamma(N-1)}.$$

\vspace{0.4 cm}

Furthermore, assuming that $N$ is sufficiently big for $\gamma(N-1)$ to be greater than one, and using again the fact that the elements of $\mu$ need to sum up to $1$, we obtain $$ \mu_{N-1} +\frac{ \mu_1 +\mu_2 +\mu_3 +...+\mu_{N-2}}{\gamma(N-1)} \geq \frac{1}{\gamma(N-1)}.$$

Combining the two previous equations we have $$\mu_{N-1} \geq \frac{1}{2\gamma(N-1)}.$$

Hence, we finally get

\begin{align} \label{eq:bound2}
\sum_{k=1}^{\left\lfloor \frac{N}{2} \right\rfloor} \mathbb{P}\left(|\widetilde{\zeta}_N (\beta_N s)| = k \right) &\leq 2\gamma(N-1) \frac{N^2}{2}\left(\frac{1+ \gamma}{\gamma}\right)^{1-\frac{N}{2}} \nonumber \\
&\leq (1+\gamma)N^3\left(\frac{\gamma}{1+\gamma}\right)^{\frac{N}{2}}. 
\end{align}

\vspace{0.4 cm}
And the last bound goes to zero as $N$ goes to infinity, so that we can find some $N_2$ such that (\ref{sumsmalln}) is satisfied for any $N \geq N_2$. Finally (\ref{sumbign}) and (\ref{sumsmalln}) are both satisfied for $N \geq \max(N_1,N_2)$ so that the proof is over.
\end{proof}

\vspace{0.4 cm}

From Proposition \ref{prop:difftozero} we obtain the following corollary.

\vspace{0.4 cm}

\begin{corollary} \label{cor:convbeta}
For any $\gamma > 0$ the following convergence holds
\begin{align}\label{betaNconv}
\frac{\sigma_N}{\beta_N} \overset{\mathcal{D}}{\underset{N \rightarrow \infty}{\longrightarrow}} \mathcal{E} (1).
\end{align}
\end{corollary}

\begin{proof}
This result follows from the definition of $\beta_N$ and a simple density argument. Taking $s = t = 1$ in Proposition \ref{prop:difftozero} we have $$ \lim_{N \rightarrow \infty} \P\left( \frac{\sigma_N}{\beta_N} > 2 \right) = e^{-2}.$$

Iterating this argument we easily obtain that for any $n \in \N$ \begin{equation} \label{convnpos}
\lim_{N \rightarrow \infty} \P\left( \frac{\sigma_N}{\beta_N} > 2^n \right) = e^{-2^n}.
\end{equation}

\vspace{0.4 cm}

Moreover, taking $s = t = 2^{-1}$ in Proposition \ref{prop:difftozero} we have $$ \lim_{N \rightarrow \infty} \P\left( \frac{\sigma_N}{\beta_N} > 2^{-1} \right) = e^{-2^{-1}}.$$

And iterating this argument once again we have for any $n \in \N$ \begin{equation} \label{convnneg}
\lim_{N \rightarrow \infty} \P\left( \frac{\sigma_N}{\beta_N} > 2^{-n} \right) = e^{-2^{-n}}.
\end{equation}

\vspace{0.4 cm}

Now for any $n \in \N$ consider the set $B_n$ of the real numbers which binary expansion contains exactly $n$ times the digit $1$, that is:
$$ B_n \mydef \left\{ x \in \R: x = \sum_{k \in \Z} a_k 2^k \text{, where } (a_k)_{k \in \Z} \in \{0,1\}^\Z \text{ is such that } \sum_{k \in \Z} a_k = n \right\}.$$

We can then reformulate equations (\ref{convnpos}) and (\ref{convnneg}) as:
$$ \lim_{N \rightarrow \infty} \P\left( \frac{\sigma_N}{\beta_N} > t \right) = e^{-t} \ \ \  \text{ for any } t \in B_1.$$

\vspace{0.4 cm}

By induction we then easily get (using Proposition \ref{prop:difftozero}) that for any $n \in \N$
$$ \lim_{N \rightarrow \infty} \P\left( \frac{\sigma_N}{\beta_N} > t \right) = e^{-t} \ \ \  \text{ for any } t \in B_n.$$

\vspace{0.4 cm}

Now notice that any real number has a binary expansion, so that it can be easily approximated—by below or by above indifferently—by a sequence in $\bigcup_{n \in \N} B_n$. Then the result follows by monotonicity of the function $t \mapsto \P\left( \frac{\sigma_N}{\beta_N} > t \right)$. 
\end{proof}

\vspace{0.4 cm}

It only remains to prove that we can replace $\beta_N$ by $\E \left( \sigma_N \right)$ in Corollary \ref{cor:convbeta}. It follows from the proposition below.

\begin{proposition} \label{prop:betaexpequi}
For any $\gamma > 0$ the expectation of $\sigma_N$ is finite and the following convergence holds

$$\lim_{N \rightarrow \infty} \frac{\E \left( \sigma_N \right)}{\beta_N} = 1.$$
\end{proposition}

\begin{proof}
From the following identity \begin{align*}
&\left[ \P\left( \frac{\sigma_N}{\beta_N} > s + t \right) - \P \left( \frac{\sigma_N}{\beta_N} > s \right)\P\left( \frac{\sigma_N}{\beta_N} > t \right) \right] \\
&= \sum_{k=1}^{N} \left[\mathbb{P}\left(\frac{\sigma_N^{\llbracket 1,k\rrbracket}}{\beta_N}>t\right)-\mathbb{P}\left(\frac{\sigma_N}{\beta_N}>t\right)\right] \cdot \mathbb{P}\big(|\zeta_N(\beta_N s)| = k \big),
\end{align*} and from the fact that for any $t \geq 0$ and $k \in \llbracket 1,N\rrbracket$ we have $$\mathbb{P}\left(\frac{\sigma_N^{\llbracket 1,k\rrbracket}}{\beta_N}>t\right) \leq \mathbb{P}\left(\frac{\sigma_N}{\beta_N}>t\right),$$ it follows that $$\P\left( \frac{\sigma_N}{\beta_N} > s + t \right) \leq \P \left( \frac{\sigma_N}{\beta_N} > s \right)\P\left( \frac{\sigma_N}{\beta_N} > t \right).$$

\vspace{0.4 cm}

Moreover from the definition of $\beta_N$ we have that for any $N \in \N$ and any $n \in \N$ 

$$ \P\left( \frac{\sigma_N}{\beta_N} > n \right) \leq e^{-n}.$$

\vspace{0.4 cm}

Therefore, for any $N \in \N$ and any $t \geq 0$ we have

$$ \P\left( \frac{\sigma_N}{\beta_N} > t \right) \leq e^{- \lfloor t \rfloor }.$$

\vspace{0.4 cm}

Firstly it implies that $\E (\sigma_N)$ is finite by a simple change of variable:
$$ \E (\sigma_N) = \int_0^{\infty} \P \left(\sigma_N > t \right) dt = \beta_N \int_0^{\infty} \P \left(\frac{\sigma_N}{\beta_N} > t \right) dt \leq \beta_N \int_0^{\infty} e^{- \lfloor t \rfloor } dt < \infty.$$

\vspace{0.4 cm}
Secondly, using the Dominated Convergence Theorem and Corollary \ref{cor:convbeta}, we get

\begin{align*}
\lim_{N \rightarrow \infty} \frac{\E \left( \sigma_N \right)}{\beta_N} &= \lim_{N \rightarrow \infty} \int_0^{\infty} \P \left(\frac{\sigma_N}{\beta_N} > t \right) dt \\
&= \int_0^{\infty} \lim_{N \rightarrow \infty} \P \left(\frac{\sigma_N}{\beta_N} > t \right) dt \\
&= \int_0^{\infty} e^{-t} dt \\
&= 1.
\end{align*}
\end{proof}

\vspace{0.4 cm}

\section*{Acknowledgements}
 This article was produced as part of the activities of FAPESP  Research, Innovation and Dissemination Center for Neuromathematics (grant number 2013/07699-0 , S.Paulo Research Foundation), and the authors were supported by  FAPESP scholarships (grant number 2017/02035-7 and 2019/14367-0).

\end{document}